\documentclass[12pt]{amsart}
\usepackage{geometry, fullpage, amssymb} 
\theoremstyle{plain}
\newtheorem{thm}{Theorem}[section]

\newtheorem{lem}[thm]{Lemma}
\newtheorem*{MainResult}{Theorem A}
\newtheorem*{lonelyPathLemma}{Lonely Path Lemma}
\newtheorem*{generalizedLonelyPathLemma}{Generalized Lonely Path Lemma}
\newtheorem{cor}[thm]{Corollary}
\newtheorem*{UnnumberedCorollary}{Corollary}
\newtheorem*{conjecture}{Conjecture}

\theoremstyle{definition}
\newtheorem{defn}{Definition}
\theoremstyle{remark}

\newtheorem*{question}{Question}

\title{Coloring and The Lonely Graph}
\author{landon rabern}

\begin{document}
\maketitle

\begin{abstract}
We improve upper bounds on the chromatic number proven independently in \cite{reedNote} and \cite{ingo}.  Our main lemma gives a sufficient condition for two paths in graph to be completely joined.  Using this, we prove that if a graph has an optimal coloring with more than  $\frac{\omega}{2}$ singleton color classes, then it satisfies $\chi \leq \frac{\omega + \Delta + 1}{2}$.  It follows that a graph satisfying $n - \Delta < \alpha + \frac{\omega - 1}{2}$ must also satisfy $\chi \leq \frac{\omega + \Delta + 1}{2}$, improving the bounds in \cite{reedNote} and \cite{ingo}.  We then give a simple argument showing that
if a graph satisfies $\chi > \frac{n + 3 - \alpha}{2}$, then it also satisfies $\chi(G) \leq \left\lceil\frac{\omega(G) + \Delta(G) + 1}{2}\right\rceil$.  From this it follows that a graph satisfying $n - \Delta < \alpha + \omega$ also satisfies $\chi(G) \leq \left\lceil\frac{\omega(G) + \Delta(G) + 1}{2}\right\rceil$ improving the bounds in \cite{reedNote} and \cite{ingo} even further at the cost of a ceiling.
In the next sections, we generalize our main lemma to constrained colorings (e.g. $r$-bounded colorings).  We present a generalization of Reed's conjecture to $r$-bounded colorings and prove the conjecture for graphs with maximal degree close to their order.  Finally, we outline some applications (in  \cite{BorodinKostochka}  and \cite{ColoringWithDoublyCriticalEdge}) of the theory presented here to the Borodin-Kostochka conjecture and coloring graphs containing a doubly critical edge.
\end{abstract}

\section{Frames and lonely edges}

The vertex swapping operation that we will study preserves the following structure of  a coloring.
\begin{defn}
Let $C = \{I_1,\ldots, I_m\}$ be a coloring of a graph $G$.  The \emph{frame} of $C$ (denoted $Frame(C)$) is the sequence $|I_{j_1}|,  |I_{j_2}|, \ldots, |I_{j_m}|$ where the $1 \leq j_k \leq m$ are distinct and $j_a \leq j_b \Rightarrow |I_{j_a}|  \leq  |I_{j_b}|$.  In other words, the ordered sequence of color class orders.  Let $|Frame(C)|$ denote the length of $Frame(C)$. Let $0$ denote the unique zero length frame.
\end{defn}

\begin{defn}
Let $C = \{I_1,\ldots, I_m\}$ be a coloring of a graph $G$.  If there exists $j \neq k$ such that $v \in I_j$, $w \in I_k$ and $N(v) \cap I_k = \{w\}$, then the (directed) edge $(v, w)$ is called \emph{$C$-lonely}.  If the coloring is clear from context we drop the $C$ and just call the edge plain \emph{lonely}.
\end{defn}

The following simple lemma is immediate from the definitions.

\begin{lem}\label{SwappingPreservesFrame}
Let $C$ be a coloring of a graph $G$.  If both $(v, w)$ and $(w,v)$ are $C$-lonely, then swapping $v$ and $w$ yields a new coloring $C'$ on the same frame.
\end{lem}

\begin{defn}
Let $C$ be a coloring of a graph $G$.  The \emph{$C$-lonely graph of $G$} (denoted $L_C(G)$) is the directed graph with vertex set $V(G)$ and edge set $\{(v, w) \mid (v, w) \text{ is $C$-lonely in $G$}\}$.
\end{defn}

The next lemma gives us a way to force dense strips in graphs with many lonely edges.  

\begin{lonelyPathLemma}
Let $G$ be a graph.  If $C$ is an optimal coloring of $G$, $\{a\}, \{b\} \in C$ are distinct singleton color classes and $p_a$, $p_b$ are vertex disjoint (directed) paths in $L_C(G)$ (starting at $a$, $b$ respectively) both having at most one vertex in any given color class, then the vertices of $p_a$ are completely joined to the vertices of $p_b$ in $G$.
\end{lonelyPathLemma}
\begin{proof}
Assume (to reach a contradiction) that the lemma is false. Of all counterexamples, pick an optimal coloring $C$ of $G$, $\{a\}, \{b\} \in C$ distinct singleton color classes and $p_a$, $p_b$ vertex disjoint (directed) paths in $L_C(G)$ (starting at $a$, $b$ respectively) both having at most one vertex in any given color class where the sum of the lengths of $p_a$ and $p_b$ is minimized.  Then, by the minimality condition, all but the ends of $p_a$ and $p_b$ must be joined in $G$.  If $p_a$ contains more than one vertex (say $p_a = a, a_2, a_3, \ldots, a_n$), then $(a, a_2)$ is lonely since $p_a$ is a path in $L_C(G)$.  But \{a\} is a singleton color class, so $(a_2, a)$ is also lonely.  Hence, by Lemma \ref{SwappingPreservesFrame}, swapping $a$ and $a_2$ yields another optimal coloring $C'$ of $G$.  

To apply the minimality condition, we need to show that $p_a' = a_2, a_3, \ldots, a_n$ and $p_b$ are paths in $L_{C'}(G)$.  Let $I_j$, $I_j'$ be the color classes containing $a_j$ in $C$, $C'$ respectively.  Assume that $p_a' \not \in L_{C'}(G)$. Then we have $2 \leq k \leq n - 1$ such that $|N(a_k) \cap I_{k+1}'| \neq 1$.  Hence $I_{k+1}' \neq I_{k+1}$. Since swapping $a$ and $a_2$ only changes $\{a\}$ and $I_2$, we must have $I_{k+1} = \{a\}$ or $I_{k+1} = I_2$.  In the latter case, $a_{k+1} = a_2$ since $p_a$ has at most one vertex in each color class. Thus $a_{k+1} = a$ or $a_{k+1} = a_2$.  If $a_{k+1} = a_2$, then $I_{k+1}' = \{a_{k+1}\}$ contradicting the fact that $|N(a_k) \cap I_{k+1}'| \neq 1$.    
Whence $a_{k+1} = a$.  Since $p_a$ is a path, it has no repeated internal vertices; hence, $k+1 = n$. This is a contradiction since $a_n$ is not joined to the end of $p_b$ but $a$ is. Whence $p_a' \in L_{C'}(G)$.

Now assume that $p_b \not \in L_{C'}(G)$ (say $p_b = b, b_2, \ldots, b_m$).  Let $Q_j$, $Q_j'$ be the color classes containing $b_j$ in $C$, $C'$ respectively. Then we have $2 \leq e \leq m - 1$ such that $|N(b_e) \cap Q_{e+1}'| \neq 1$. Hence $Q_{e+1}' \neq Q_{e+1}$. Since swapping $a$ and $a_2$ only changes $\{a\}$ and $I_2$, we must have $Q_{e+1} = \{a\}$ or $Q_{e+1} = I_2$.  The former is impossible since $p_a$ and $p_b$ are disjoint.  Hence $Q_{e+1} = I_2$.  Since $e < m$, $b_e$ is adjacent to $a_2$.  Since $|N(b_e) \cap I_2| = |N(b_e) \cap Q_{e+1}| = 1$, we must have $b_{e+1} = a_2$ contradicting the disjointness of $p_a$ and $p_b$. Whence $p_b \in L_{C'}(G)$.

Hence $p_a'$ and $p_b$ are vertex disjoint paths in $L_{C'}(G)$ with the end of $p_a'$ not joined to the end of $p_b$ and $p_a'$ shorter than $p_a$, contradicting the minimality condition.  Hence $p_a$ is the single vertex $\{a\}$.  Similarly, $p_b$ is the single vertex $\{b\}$.  Since $p_a$ is not joined to $p_b$, the color classes $\{a\}$ and $\{b\}$ can be merged, contradicting the fact that $C$ is an optimal coloring. 
\end{proof}

The end of this section shows that graphs which do not satisfy Reed's $\omega$, $\Delta$, and $\chi$ bound are replete with lonely edges.

\begin{defn}
Let $C$ be a coloring of a graph $G$.  For any vertex $v \in G$, set
\[L_C(v) = \{w \in G \mid (v,w) \text{ is C-lonely}\}.\] 
\end{defn}

\begin{lem}\label{SomebodyTouchesEverybody}
Let $G$ be a graph and $C = \{I_1,\ldots,I_m\}$ an optimal coloring of $G$. Then, for each $1 \leq j \leq m$, there exists $v_j \in I_j$ such that $N(v_j) \cap I_k \neq \emptyset$ for each $k \neq j$.
\end{lem}
\begin{proof}
Otherwise $C$ would not be optimal.
\end{proof}

\begin{lem}\label{ReedCounterexampleRepleteWithLonelyEdges}
Let $G$ be a graph with $\chi(G) > \frac{\omega(G) + \Delta(G) + 1}{2} + t$ and $C = \{I_1,\ldots,I_m\}$ 
an optimal coloring of $G$.  Then, for each $1 \leq j \leq m$, there exists $v_j \in I_j$ such that
\[|L_C(v_j)| \geq \omega(G) + 2t.\]
\end{lem}
\begin{proof}
Fix $j$ with $1 \leq j \leq m$.  By Lemma \ref{SomebodyTouchesEverybody}, we have $v_j \in I_j$ such that $|N(v_j) \cap I_k| \geq 1$ for each $k \neq j$. Hence
\begin{align*}
d(v_j) & \geq 2(m - 1 - |L_C(v_j)|) + |L_C(v_j)| \\
&= 2m - |L_C(v_j)| - 2 \\
&= 2\chi(G) - |L_C(v_j)| - 2 
\end{align*}
But $\chi(G) > \frac{\omega(G) + \Delta(G) + 1}{2} + t$, thus
\begin{align*}
2\chi(G) - |L_C(v_j)| - 2 &\leq d(v_j) \\
&\leq \Delta(G) \\
&< 2\chi(G) - \omega(G) - (2t + 1).
\end{align*}
The lemma follows.
\end{proof}

\section{Very stingy graphs}

\begin{defn}
The \emph{stinginess} of a graph $G$ (denoted $\iota(G)$) is the maximum number of singleton color classes appearing in an optimal coloring of $G$.  An optimal coloring of $G$ is called \emph{stingy} just in case it has the maximum number of singleton color classes.
\end{defn}

\begin{lem}\label{VeryStingyImpliesReed}
If $G$ is a graph with $\iota(G) > \frac{\omega(G)}{2}$, then
\[\chi(G) \leq \frac{\omega(G) + \Delta(G) + 1}{2}.\]
\end{lem}
\begin{proof}
Assume (to reach a contradiction) that the lemma is false and let $G$ be a counterexample.  Let $C$ be a stingy coloring of $G$ and let $S$ be the vertices in singleton color classes of $C$.  By Lemma \ref{ReedCounterexampleRepleteWithLonelyEdges}, $|L_C(v)| \geq \omega(G)$ for each $v \in S$.  If there exists $v\in S$ such that $L_C(v) \subseteq \displaystyle\bigcup_{w \in S \smallsetminus \{v\}} L_C(w)$, then $\{v\} \cup L_C(v)$ induces a clique in $G$ by the Lonely Path Lemma.  But $|\{v\} \cup L_C(v)| \geq \omega(G) + 1$, so this is impossible.  Hence, for each $v \in S$, we have $l_v \in L_C(v)$ such that $l_v \not \in \displaystyle\bigcup_{w \in S \smallsetminus \{v\}} L_C(w)$.  Set $T = S \cup \displaystyle\bigcup_{v \in S} l_v$.  Then $T$ induces a clique in $G$ by the Lonely Path Lemma.  But $S \cap \displaystyle\bigcup_{v \in S} l_v = \emptyset$ and thus 
\[|T| = |S \cup \displaystyle\bigcup_{v \in S} l_v| = |S| + |\displaystyle\bigcup_{v \in S} l_v| = 2|S| = 2\iota(G) > 2\frac{\omega(G)}{2} = \omega(G).\]
This contradiction completes the proof.
\end{proof}

Note that our application of the Lonely Path Lemma was restricted to paths of length at most one.  We think that it is possible to prove better results along these lines by using the full power of the lemma.

\section{Improvements to the graph associations bound}

\begin{lem}\label{StinginessPreservedInPatchedColorings}
Let $G$ be a graph and $H$ an induced subgraph of $G$.  If $\chi(G) = \chi(G \smallsetminus H) + \chi(H)$, then $\iota(G) \geq  \iota(G \smallsetminus H) + \iota(H)$.
\end{lem}
\begin{proof}
Assume that $\chi(G) = \chi(G \smallsetminus H) + \chi(H)$.  Then patching together any optimal coloring of $G \smallsetminus H$ with any optimal coloring of $H$ yields an optimal coloring of $G$.  The lemma follows.
\end{proof}

\begin{lem}\label{ChiAtMostAverageOfStinginessAndOrder}
Let $G$ be a graph. Then $\chi(G) \leq \frac{\iota(G) + |G|}{2}$.
\end{lem}
\begin{proof}
Let $C = \{I_1,\ldots,I_m, \{s_1\}, \ldots, \{s_{\iota(G)}\}\}$ be a stingy coloring of $G$.  Since $|I_j| \geq 2$ for $1 \leq j \leq m$, we have $\chi(G) \leq \iota(G) + \frac{|G| - \iota(G)}{2} = \frac{|G| + \iota(G)}{2}$.
\end{proof}

In \cite{graphAssociations} and \cite{ingo}, the bound $\chi(G) \leq \frac{\omega(G) + |G| - \alpha(G) + 1}{2}$ was proven.  The following improves this bound.

\begin{thm}\label{ReedDisjunct}
For any graph $G$, at least one of the following holds,
\begin{enumerate}
\item $\chi(G) \leq \frac{\omega(G) + \Delta(G) + 1}{2}$,
\item $\chi(G) \leq \frac{\frac{\omega(G)}{2} + |G| - \alpha(G)}{2} + 1$.
\end{enumerate}
\end{thm}
\begin{proof}
Assume (to reach a contradiction) that this is not the case and let $G$ be a counterexample with the minimum number of vertices.
Let $I$ be a maximum independent set in $G$. Then $\chi(G \smallsetminus I) \leq \chi(G) \leq \chi(G \smallsetminus I) + 1$.  Since $|G \smallsetminus I| < |G|$, the theorem holds for $G \smallsetminus I$. Hence $\chi(G) = \chi(G \smallsetminus I) + 1$.  Whence, by Lemma \ref{StinginessPreservedInPatchedColorings}, we have $\iota(G) \geq \iota(G \smallsetminus I) $.  Assume that (1) does not hold for $G$.  Then, by Lemma \ref{VeryStingyImpliesReed}, $\iota(G \smallsetminus I) \leq \iota(G) \leq \frac{\omega(G)}{2}$. Hence, by Lemma \ref{ChiAtMostAverageOfStinginessAndOrder}, we have 
\begin{align*}
\chi(G) &= \chi(G \smallsetminus I) + 1 \\
&\leq \frac{|G| - |I| + \iota(G)}{2} + 1 \\
&\leq \frac{|G| - \alpha(G) + \frac{\omega(G)}{2}}{2} + 1.
\end{align*}
\end{proof}

In both \cite{reedNote} and \cite{ingo} it was proven  that if $\chi(G) > \frac{\omega(G) + \Delta(G) + 1}{2}$, then $|G| - \Delta(G) \geq \alpha(G) + 1$. Using Theorem \ref{ReedDisjunct}, we can easily deduce an improvement of this bound.

\begin{cor}
If $G$ is a graph satisfying $\chi(G) > \frac{\omega(G) + \Delta(G) + 1}{2}$, then 
\[|G| - \Delta(G) \geq \alpha(G) + \frac{\omega(G) - 1}{2}.\]
\end{cor}
\begin{proof}
Let $G$ be such a graph.  By Theorem \ref{ReedDisjunct}, 
\[ \frac{\omega(G) + \Delta(G) + 1}{2} < \chi(G) \leq \frac{\frac{\omega(G)}{2} + |G| - \alpha(G)}{2} + 1.\]
Hence $|G| - \Delta(G) > \alpha(G) + \frac{\omega(G)}{2} - 1$.  The corollary follows.
\end{proof}

\section{A Cheap Improvement}
The following two lemmas were proved in \cite{reedNote} using matching theory results.
\begin{lem}\label{ReedForChiBIggerThanHalf}
If $G$ is a graph with $\chi(G) > \left \lceil \frac{|G|}{2} \right \rceil$, then 
\[\chi(G) \leq \frac{\omega(G) + \Delta(G) + 1}{2}.\]
\end{lem}

\begin{lem}\label{ReedForAlphaAtMostTwo}
If $G$ is a graph with $\alpha(G) \leq 2$, then
\[\chi(G) \leq \left\lceil\frac{\omega(G) + \Delta(G) + 1}{2}\right\rceil.\]
\end{lem}

The following simple bound is proved by just pulling out a maximal independent set and seeing what happens.
\begin{thm}\label{SimpleBound}
If $G$ is a graph with $\chi(G) > \frac{|G| + 3 - \alpha(G)}{2}$, then
\[\chi(G) \leq \left\lceil\frac{\omega(G) + \Delta(G) + 1}{2}\right\rceil.\]
\end{thm}
\begin{proof}
Let $G$ be a graph with $\chi(G) > \frac{|G| + 3 - \alpha(G)}{2}$ and $I$ an independent set in $G$ with $\alpha(G)$ vertices.  Put $H = G \smallsetminus I$.  Then 
\begin{align*}
\chi(H) &\geq \chi(G) - 1 \\
&> \frac{|G| + 3 - \alpha(G)}{2} - 1 \\
&= \frac{|G| + 1 - \alpha(G)}{2} \\
&=  \frac{|H| + 1}{2}.
\end{align*}

Hence, by Lemma \ref{ReedForChiBIggerThanHalf}, we have 
\[\chi(H) \leq \frac{\omega(H) + \Delta(H) + 1}{2}.\]

But $I$ is a maximal independent set and hence each vertex of $H$ is adjacent to at least one vertex in $I$.  In particular, $\Delta(H) \leq \Delta(G) - 1$.  Whence

\[\chi(G) \leq \chi(H) + 1 \leq \frac{\omega(H) + \Delta(G) - 1 + 1}{2} + 1 \leq  \frac{\omega(G) + \Delta(G) + 1}{2} + \frac{1}{2}.\]
The theorem follows.
\end{proof}

\begin{cor}
If $G$ is a graph with $\chi(G) \geq \frac{|G| + 1}{2}$, then 
\[\chi(G) \leq \left\lceil\frac{\omega(G) + \Delta(G) + 1}{2}\right\rceil.\]
\end{cor}
\begin{proof}
Let $G$ be a graph with $\chi(G) \geq \frac{|G| + 1}{2}$.  If $\alpha(G) \leq 2$, then we are done by Lemma \ref{ReedForAlphaAtMostTwo}. If $\alpha(G) \geq 3$, then $\chi(G) \geq \frac{|G| + 1}{2} > \frac{|G|}{2} \geq \frac{|G| + 3 - \alpha(G)}{2}$ and we are done by Theorem \ref{SimpleBound}. 
\end{proof}

\begin{cor}
If $G$ is a graph with $\chi(G) >  \left\lceil\frac{\omega(G) + \Delta(G) + 1}{2}\right\rceil$, then 
\[|G| - \Delta(G) \geq \alpha(G) +\omega(G).\]
\end{cor}
\begin{proof}
Let $G$ be such a graph.  By Theorem \ref{SimpleBound}, 
\[ \frac{\omega(G) + \Delta(G) + 1}{2} + \frac{1}{2} < \chi(G) \leq \frac{|G| + 3 - \alpha(G)}{2}.\]
Hence $|G| - \Delta(G) > \alpha(G) + \omega(G) - 1$.  The corollary follows.
\end{proof}

\section{A generalization of the Lonely Path Lemma}
\begin{defn}
Let $G$ be a graph. A \emph{property of colorings on $G$} is a subset of the set of all (proper) colorings of $G$.
\end{defn}

\begin{defn}
A property $P$ of colorings on a graph $G$ is a \emph{frame property} just in case 
\[Frame(C) = Frame(C') \Rightarrow [C \in P \Rightarrow C' \in P],\]
for any colorings $C$, $C'$ of $G$.
\end{defn}

\begin{defn}
A property $P$ of colorings on a graph $G$ is \emph{singleton-friendly} just in case
\[C \in P \Rightarrow C' \in P,\]
 for any coloring $C'$ formed by merging two singleton color classes of a coloring $C$. 
\end{defn}

\begin{defn}
Let $P$ be a property of colorings on a graph $G$.  A coloring $C$ of $G$ is \emph{$P$-optimal}  just in case $|C|$ is minimal among colorings of $G$ satisfying $P$.  Let $\chi_P(G)$ denote the order of a $P$-optimal coloring of $G$.
\end{defn}

\begin{generalizedLonelyPathLemma}
Let $G$ be a graph and $P$ a singleton-friendly frame property.  If $C$ is a $P$-optimal coloring of $G$, $\{a\}, \{b\} \in C$ are distinct singleton color classes and $p_a$, $p_b$ are vertex disjoint (directed) paths in $L_C(G)$ (starting at $a$, $b$ respectively) both having at most one vertex in any given color class, then the vertices of $p_a$ are completely joined to the vertices of $p_b$ in $G$.
\end{generalizedLonelyPathLemma}
\begin{proof}
Assume (to reach a contradiction) that the lemma is false. Of all counterexamples, pick a $P$-optimal coloring $C$ of $G$, $\{a\}, \{b\} \in C$ distinct singleton color classes and $p_a$, $p_b$ vertex disjoint (directed) paths in $L_C(G)$ (starting at $a$, $b$ respectively) both having at most one vertex in any given color class where the sum of the lengths of $p_a$ and $p_b$ is minimized.  Then, by the minimality condition, all but the ends of $p_a$ and $p_b$ must be joined in $G$.  If $p_a$ contains more than one vertex (say $p_a = a, a_2, a_3, \ldots$), then $(a, a_2)$ is lonely since $p_a$ is a path in $L_C(G)$.  But \{a\} is a singleton color class, so $(a_2, a)$ is also lonely.  Hence, by Lemma \ref{SwappingPreservesFrame}, swapping $a$ and $a_2$ yields a new coloring $C'$ on the same frame.  Since $P$ is a frame property, $C'$ is $P$-optimal.  

To apply the minimality condition, we need to show that $p_a' = a_2, a_3, \ldots, a_n$ and $p_b$ are paths in $L_{C'}(G)$.  Let $I_j$, $I_j'$ be the color classes containing $a_j$ in $C$, $C'$ respectively.  Assume that $p_a' \not \in L_{C'}(G)$. Then we have $2 \leq k \leq n - 1$ such that $|N(a_k) \cap I_{k+1}'| \neq 1$.  Hence $I_{k+1}' \neq I_{k+1}$. Since swapping $a$ and $a_2$ only changes $\{a\}$ and $I_2$, we must have $I_{k+1} = \{a\}$ or $I_{k+1} = I_2$.  In the latter case, $a_{k+1} = a_2$ since $p_a$ has at most one vertex in each color class. Thus $a_{k+1} = a$ or $a_{k+1} = a_2$.  If $a_{k+1} = a_2$, then $I_{k+1}' = \{a_{k+1}\}$ contradicting the fact that $|N(a_k) \cap I_{k+1}'| \neq 1$.    
Whence $a_{k+1} = a$.  Since $p_a$ is a path, it has no repeated internal vertices; hence, $k+1 = n$. This is a contradiction since $a_n$ is not joined to the end of $p_b$ but $a$ is. Whence $p_a' \in L_{C'}(G)$.

Now assume that $p_b \not \in L_{C'}(G)$ (say $p_b = b, b_2, \ldots, b_m$).  Let $Q_j$, $Q_j'$ be the color classes containing $b_j$ in $C$, $C'$ respectively. Then we have $2 \leq e \leq m - 1$ such that $|N(b_e) \cap Q_{e+1}'| \neq 1$. Hence $Q_{e+1}' \neq Q_{e+1}$. Since swapping $a$ and $a_2$ only changes $\{a\}$ and $I_2$, we must have $Q_{e+1} = \{a\}$ or $Q_{e+1} = I_2$.  The former is impossible since $p_a$ and $p_b$ are disjoint.  Hence $Q_{e+1} = I_2$.  Since $e < m$, $b_e$ is adjacent to $a_2$.  Since $|N(b_e) \cap I_2| = |N(b_e) \cap Q_{e+1}| = 1$, we must have $b_{e+1} = a_2$ contradicting the disjointness of $p_a$ and $p_b$. Whence $p_b \in L_{C'}(G)$.

Thence $p_a'$ and $p_b$ are vertex disjoint paths in $L_{C'}(G)$ with the end of $p_a'$ not joined to the end of $p_b$ and $p_a'$ shorter than $p_a$, contradicting the minimality condition.  Hence $p_a$ is the single vertex $\{a\}$.  Similarly, $p_b$ is the single vertex $\{b\}$.  Since $p_a$ is not joined to $p_b$, the color classes $\{a\}$ and $\{b\}$ can be merged to yield a new coloring $D$.  Since $P$ is singleton-friendly, $D$ satisfies $P$.  But $|D| < |C|$, contradicting the fact that $C$ is $P$-optimal.
\end{proof}

Before we can do anything with this lemma, we need to find some interesting singleton-friendly frame properties.

\begin{question}
What does a singleton-friendly frame property look like?
\end{question}

There is a simple sufficient condition for a property to be a singleton-friendly frame property.

\begin{defn}
Let $C$ be a coloring.  Denote by $Frame_m(C)$ the subsequence of $Frame(C)$ beginning with the first $m$.
\end{defn}

\begin{lem}\label{SingletonFriendlyFramePropertySufficientCondition}
Let $P$ be a property of colorings on a graph $G$.  If 
\[Frame_3(C) = Frame_3(C')  \Rightarrow [C \in P \Rightarrow C' \in P],\]
for any colorings $C$, $C'$ of $G$, then $P$ is a singleton-friendly frame property.
\end{lem}
\begin{proof}
Assume that $Frame_3(C) = Frame_3(C')  \Rightarrow [C \in P \Rightarrow C' \in P]$ for any colorings $C$, $C'$ of $G$.  Plainly, $P$ is a frame property.  Since merging singleton color classes only affects the $1$'s and $2$'s of a frame, we see that $P$ is also singleton-friendly.
\end{proof}

This condition is not necessary.  For example, consider the property ``has at most $k$ singleton color classes''.  The condition can be made sufficient by considering the total number of vertices in singleton and doubleton color classes.

\begin{defn}
Given a coloring $C$ of a graph $G$, let $Small(C)$ be the order of the union of the singleton and doubleton color classes of $C$.
\end{defn}

\begin{lem}\label{SingletonFriendlyFramePropertyCompleteCondition}
Let $P$ be a property of colorings on a graph $G$.  Then $P$ is a singleton-friendly frame property if and only if
\[[Small(C) = Small(C') \wedge Frame_3(C) = Frame_3(C')]  \Rightarrow [C \in P \Rightarrow C' \in P],\]
for any colorings $C$, $C'$ of $G$.
\end{lem}

The following two lemmas, which are immediate from the definitions, describe the basic structure of the properties under consideration.
\begin{lem}
Let $G$ be a graph. The frame properties on $G$ are a topology on the set of (proper) colorings of $G$.
\end{lem}

\begin{lem}
Let $G$ be a graph. The singleton-friendly frame properties on $G$ are a topology on the set of (proper) colorings of $G$.
\end{lem}

\section{Reed's conjecture generalized to $r$-bounded colorings}
\begin{defn}
Let $G$ be a graph and $r$ a natural number. An \emph{$r$-bounded coloring} of $G$ is a (proper) coloring of $G$ in which all color classes have order at most $r$.  
\end{defn}

Observe that a coloring $C$ is an $r$-bounded coloring of a graph $G$ just in case $Frame_{r+1}(C) = 0$.

\begin{lem}
Let $G$ be a graph and $r \geq 2$. Let $B_r = \{C \mid \text{$C$ is an $r$-bounded coloring of $G$}\}$.  Then $B_r$ is a singleton-friendly frame property.
\end{lem}
\begin{proof}
Let $C \in B_r$ and $C'$ be a coloring of $G$ with $Frame_3(C) = Frame_3(C')$.  Then, since $r + 1 \geq 3$, $Frame_{r+1}(C') = Frame_{r+1}(C)$. Also, since $C$ is $r$-bounded,
\[Frame_{r+1}(C') = Frame_{r+1}(C) = 0.\]

Thus $C'$ is $r$-bounded as well and we have $C' \in B_r$.
Hence the lemma follows from Lemma \ref{SingletonFriendlyFramePropertySufficientCondition}.
\end{proof}

To simplify notation a bit, we write $\chi_r(G)$ in place of $\chi_{B_r}(G)$.

\begin{lem}\label{SingletonsTouchEverybodyExceptMaybeTheBigBoys}
Let $C = \{I_1,\ldots,I_m\}$ be an optimal $r$-bounded coloring of a graph $G$.  If $I_j = \{v\}$ for some $j$, then $N(v) \cap I_k \neq \emptyset$ for each $k \neq j$ such that $|I_k|< r$.
\end{lem}
\begin{proof}
Otherwise $C$ would not be optimal.
\end{proof}

\begin{defn}
Let $G$ be a graph.  Denote the maximum number of order $r$ color classes in an optimal $r$-bounded coloring of $G$ by $M_r(G)$.  That is,

\[M_r(G) = \max\{|Frame_r(C)| \mid \text{C is an optimal $r$-bounded coloring of $G$}\}.\]
\end{defn}

\begin{lem}\label{GeneralizedReedCounterexampleRepleteWithLonelyEdges}
Let $G$ be a graph with $\chi_r(G) - M_r(G) > \frac{\omega(G) + \Delta(G) + 1}{2} + t$ and $C = \{I_1,\ldots,I_m\}$ 
an optimal $r$-bounded coloring of $G$.  If $I_j = \{v\}$ for some $j$, then
\[|L_C(v)| \geq \omega(G) + 2t.\]
\end{lem}
\begin{proof}
Assume that $I_j = \{v\}$. By Lemma \ref{SingletonsTouchEverybodyExceptMaybeTheBigBoys}, $|N(v) \cap I_k| \geq 1$ for each $k \neq j$ such that $|I_k|< r$. There are precisely $m - |Frame_r(C)| - 1$ such $k$; hence
\begin{align*}
d(v) & \geq 2(m - |Frame_r(C)| - 1 - |L_C(v)|) + |L_C(v)| \\
&= 2(m-|Frame_r(C)|) - |L_C(v)| - 2 \\
&= 2(\chi_r(G) - |Frame_r(C)|) - |L_C(v)| - 2 
\end{align*}
But $\chi_r(G) - M_r(G) > \frac{\omega(G) + \Delta(G) + 1}{2} + t$, thus
\begin{align*}
2(\chi_r(G) - M_r(G)) - |L_C(v)| - 2 &\leq 2(\chi_r(G) - |Frame_r(C)|) - |L_C(v)| - 2 \\
&\leq d(v) \\
&\leq \Delta(G) \\
&< 2(\chi_r(G) -M_r(G)) - \omega(G) - (2t + 1).
\end{align*}
The lemma follows.
\end{proof}

\begin{defn}
The \emph{$r$-bounded stinginess} of a graph $G$ (denoted $\iota_r(G)$) is the maximum number of singleton color classes appearing in an optimal $r$-bounded coloring of $G$.  An optimal $r$-bounded coloring of $G$ is called \emph{stingy} just in case it has the maximum number of singleton color classes.
\end{defn}

\begin{thm}\label{VeryStingyImpliesGeneralizedReed}
If $G$ is a graph with $\iota_r(G) > \frac{\omega(G)}{2}$, then
\[\chi_r(G) - M_r(G) \leq \frac{\omega(G) + \Delta(G) + 1}{2}.\]
\end{thm}
\begin{proof}
Assume (to reach a contradiction) that the lemma is false and let $G$ be a counterexample.  Let $C$ be a stingy $r$-bounded coloring of $G$ and let $S$ be the vertices in singleton color classes of $C$.  By Lemma \ref{GeneralizedReedCounterexampleRepleteWithLonelyEdges}, $|L_C(v)| \geq \omega(G)$ for each $v \in S$.  If there exists $v\in S$ such that $L_C(v) \subseteq \displaystyle\bigcup_{w \in S \smallsetminus \{v\}} L_C(w)$, then $\{v\} \cup L_C(v)$ induces a clique in $G$ by the Generalized Lonely Path Lemma.  But $|\{v\} \cup L_C(v)| \geq \omega(G) + 1$, so this is impossible.  Hence, for each $v \in S$, we have $l_v \in L_C(v)$ such that $l_v \not \in \displaystyle\bigcup_{w \in S \smallsetminus \{v\}} L_C(w)$.  Set $T = S \cup \displaystyle\bigcup_{v \in S} l_v$.  Then $T$ induces a clique in $G$ by the Generalized Lonely Path Lemma.  But $S \cap \displaystyle\bigcup_{v \in S} l_v = \emptyset$ and thus 
\[|T| = |S \cup \displaystyle\bigcup_{v \in S} l_v| = |S| + |\displaystyle\bigcup_{v \in S} l_v| = 2|S| = 2\iota_r(G) > 2\frac{\omega(G)}{2} = \omega(G).\]
This contradiction completes the proof.
\end{proof}

Since $\chi_2(G) - M_2(G) = \iota_2(G)$ we can drop the $\iota_r(G) > \frac{\omega(G)}{2}$ condition for the $r=2$ case and conclude the following.
\begin{cor}\label{iota_2Works}
If $G$ is a graph, then
\[\iota_2(G) \leq \frac{\omega(G) + \Delta(G) + 1}{2}.\]
\end{cor}

We rewrite this corollary in terms of standard graph properties. 
\begin{defn}
The \emph{matching number} of a graph $G$, denoted $\nu(G)$ is the maximum number of edges in a matching of $G$.
\end{defn}

Note that $\iota_2(G) = |G| - 2\nu(\overline{G})$.

\begin{cor}
If $G$ is a graph, then
\[\nu(G) \geq \frac{|G| - \alpha(G) + \delta(G)}{4}.\]
\end{cor}
\begin{proof}
Apply Corollary \ref{iota_2Works} to $\overline{G}$ to get
\begin{align*} 
|G| - 2\nu(G) = \iota_2(\overline{G}) &\leq \frac{\omega(\overline{G}) + \Delta(\overline{G}) + 1}{2} \\
&= \frac{\alpha(G) + |G| - \delta(G)}{2}.\\
\end{align*}
The corollary follows.
\end{proof}

\begin{conjecture}\label{GeneralizedReed}
If $G$ is a graph and $r$ is a natural number, then
\[\chi_r(G) - M_r(G) \leq \left\lceil\frac{\omega(G) + \Delta(G) + 1}{2}\right\rceil.\]
\end{conjecture}

This holds (trivially) for $r = 1$ since $\chi_1(G) = M_1(G) = |G|$.  By Corollary \ref{iota_2Works}, the conjecture also holds for $r=2$.  The case $r = \alpha(G) + 1$ is Reed's conjecture.\newline

In support of this conjecture, we prove it for graphs having maximal degree close to their order.

\begin{lem}\label{r-BoundedStinginessPreservedInPatchedColorings}
Let $G$ be a graph and $H$ an induced subgraph of $G$.  If $\chi_r(G) = \chi_r(G \smallsetminus H) + \chi_r(H)$, then $\iota_r(G) \geq \iota_r(G \smallsetminus H) + \iota_r(H)$.
\end{lem}
\begin{proof}
Assume that $\chi_r(G) = \chi_r(G \smallsetminus H) + \chi_r(H)$.  Then patching together any optimal $r$-bounded coloring of $G \smallsetminus H$ with any optimal $r$-bounded coloring of $H$ yields an optimal coloring of $G$.  The lemma follows.
\end{proof}

\begin{lem}\label{Chi_rAtMostAverageOfStinginessAndOrder}
Let $G$ be a graph. Then $\chi_r(G) \leq \frac{\iota_r(G) + |G|}{2}$.
\end{lem}
\begin{proof}
Let $C = \{I_1,\ldots,I_m, \{s_1\}, \ldots, \{s_{\iota_r(G)}\}\}$ be a stingy $r$-bounded coloring of $G$.  Since $|I_j| \geq 2$ for $1 \leq j \leq m$, we have $\chi_r(G) \leq \iota_r(G) + \frac{|G| - \iota_r(G)}{2} = \frac{|G| + \iota_r(G)}{2}$.
\end{proof}

\begin{thm}\label{GeneralizedReedDisjunct}
Let $G$ be a graph. Then at least one of the following holds,
\begin{enumerate}
\item $\chi_r(G) - M_r(G) \leq \frac{\omega(G) + \Delta(G) + 1}{2}$,
\item $\chi_r(G) - M_r(G) \leq \frac{\frac{\omega(G)}{2} + |G| - rM_r(G)}{2}$.
\end{enumerate}
\end{thm}
\begin{proof}
Let $C$ be an optimal $r$-bounded coloring of $G$ with $M_r(G)$ color classes of order $r$ (say $I_1, \ldots, I_{M_r(G)}$).  Set $H = \cup I_j$.  Then $\chi_r(G) = \chi_r(G \smallsetminus H) + \chi_r(H)$.

 Whence, by Lemma \ref{r-BoundedStinginessPreservedInPatchedColorings}, we have $\iota_r(G) \geq \iota_r(G \smallsetminus H) $.  Assume that (1) does not hold for $G$.  Then, by Lemma \ref{VeryStingyImpliesGeneralizedReed}, $\iota_r(G \smallsetminus H) \leq \iota_r(G) \leq \frac{\omega(G)}{2}$. Hence, by Lemma \ref{Chi_rAtMostAverageOfStinginessAndOrder}, we have 
\begin{align*}
\chi_r(G) &= \chi_r(G \smallsetminus H) + \chi_r(H) \\
&\leq \frac{|G| - |H| + \iota_r(G)}{2} + \chi_r(H) \\
&\leq \frac{|G| - rM_r(G) + \frac{\omega(G)}{2}}{2} + \chi_r(H) \\
&= \frac{|G| - rM_r(G) + \frac{\omega(G)}{2}}{2} + M_r(G).
\end{align*}
\end{proof}

\begin{cor}
If $G$ is a graph satisfying $\chi_r(G) - M_r(G) > \frac{\omega(G) + \Delta(G) + 1}{2}$, then 
\[|G| - \Delta(G) \geq rM_r(G) + \frac{\omega(G) + 3}{2}.\]
\end{cor}
\begin{proof}
Let $G$ be such a graph.  By Theorem \ref{GeneralizedReedDisjunct}, 
\[ \frac{\omega(G) + \Delta(G) + 1}{2} < \chi_r(G) - M_r(G) \leq \frac{\frac{\omega(G)}{2} + |G| - rM_r(G)}{2}.\]
Hence $|G| - \Delta(G) > rM_r(G) + \frac{\omega(G)}{2} + 1$.  The corollary follows.
\end{proof}

\section{Applications}
\begin{defn}
Let $G$ be a graph.  An edge $ab \in G$ is \emph{doubly critical} just in case $\chi(G \smallsetminus \{a, b\}) = \chi(G) - 2$.
\end{defn}

Note that a graph has a doubly critical edge if and only if $\iota(G) \geq 2$.  In \cite{BorodinKostochka} the following is proved using the Lonely Path Lemma.

\begin{MainResult}
Let $G$ be a graph containing a doubly critical edge.  If $G$ satisfies $\chi \geq \Delta \geq 9$, then $G$ contains a $K_{\Delta}$.
\end{MainResult}

This settles the following conjecture of Borodin and Kostochka for graphs containing a doubly critical edge.

\begin{conjecture}
Every graph satisfying $\chi \geq \Delta \geq  9$ contains a $K_{\Delta}$.
\end{conjecture}

Here are a couple of interesting corollaries from \cite{ColoringWithDoublyCriticalEdge}.

\begin{UnnumberedCorollary}
Let $G$ be a claw-free graph containing a doubly critical edge.  Then
\[\chi(G) \leq \left\lceil\frac{\omega(G) + \Delta(G) + 1}{2}\right\rceil.\]
\end{UnnumberedCorollary}

\begin{UnnumberedCorollary}
Let $G$ be a graph containing a doubly critical edge.  Then
\[\chi(G) \leq \textstyle\frac{1}{3}\omega(G) + \textstyle\frac{2}{3}(\Delta(G) + 1).\]
\end{UnnumberedCorollary}

\begin{question}
What does a graph containing a doubly critical edge look like?
\end{question}

\end{document}